\documentclass[12pt]{amsart}
\usepackage{amscd}

\newtheorem{thm}{Theorem}[section]
\newtheorem{cor}[thm]{Corollary}
\newtheorem{lem}[thm]{Lemma}
\newtheorem{prop}[thm]{Proposition}
\theoremstyle{definition}
\newtheorem{defn}[thm]{Definition}

\newtheorem{rem}[thm]{Remark}

\theoremstyle{remark}

\numberwithin{equation}{section}

\newcommand{\norm}[1]{\lVert#1\rVert}
\newcommand{\GL}{\mathit{GL}}

\DeclareMathOperator{\tr}{tr}

\def\bs#1{\boldsymbol{#1}}
\def\HLINE{\noalign{\hrule height1pt}}
\newcommand{\inner}[2]{\left\langle{#1},{#2}\right\rangle}
\makeatletter
\def\varin{\mathrel{\mathpalette\@varin\relax}}
\def\@varin#1{%
   \hbox{\setbox\z@\hbox{\m@th$#1\cup$}%
       \def\reserved@a{bold}%
       \dimen@\ifx\reserved@a\math@version .3\else .2\fi\p@
       \kern.5\wd\z@\kern-\dimen@
       \vrule\@width2\dimen@\@height1.08\ht\z@\@depth\z@
       \kern-\dimen@\kern-.5\wd\z@
       \box\z@}}
\makeatother

\dedicatory{Dedicated to
the memory of Jean Leray (1906-1998)}

\begin{document}
                               
\title[Cayley projective plane]
{A K\"ahler structure on the punctured cotangent bundle of the Cayley 
projective plane}





\author{K. Furutani}
\address{Department of Mathematics\\
Science University of Tokyo\\
Noda, Chiba 278-8510 Japan}
\email{furutani@ma.noda.sut.ac.jp}


\begin{abstract}
We construct a K\"ahler structure on the punctured cotangent bundle of
the Cayley projective plane whose K\"ahler form coincides with the natural
symplectic form on the cotangent bundle and we show that the geodesic
flow action is holomorphic and is expressed in a quite explicit form. 
We also give an embedding of the punctured cotangent
bundle of the Cayley projective plane into the space of $8\times 8$
complex matrices. 
\end{abstract}

\maketitle

\tableofcontents

\bigskip

\section*{Introduction}


In the paper  \cite{Ra1} a K\"ahler structure on the punctured cotangent bundle
$T^*_0S^n$ = $T^*S^n\backslash S^n$ of the sphere $S^n$ is
constructed through the mapping $\tau_S$
\begin{equation}
\begin{array}{llll}
\tau_S: & T_0^* S^n& \longrightarrow &\mathbb{C}^{n+1}\\
 & (x,y) & \longmapsto & \norm{y}x+\sqrt{-1}y.
\end{array}
\end{equation}
It is shown that the natural symplectic form $\omega$ = $\omega_{S}$ on 
the cotangent bundle coincides with the K\"ahler form ${\sqrt{-2}}\,
\overline{\partial}\partial \norm z $(see also \cite{So}).  
Moreover the geodesic flow action is
holomorphic.

In the paper  \cite{FT} we constructed a  K\"ahler structure
on the punctured cotangent bundle of complex and quaternion
projective spaces with similar properties 
as for the sphere cases(see also \cite{Ii2}).
This K\"ahler structure is just a positive complex 
polarization on the cotangent bundle as a symplectic manifold and
is applied to construct a quantization
operator by the method of pairing of polarizations.  The operator
quantizes geodesic flows of such manifolds. In other words such an
operator gives a correspondence between the geodesic flow and the one
parameter group of Fourier integral operators generated by the
square root of the Laplacian(\cite{FY}, \cite{Ra2}).

In this paper we construct a K\"ahler structure on the punctured 
cotangent bundle of the
Cayley projective plane whose K\"ahler form coincides with the natural
symplectic form on the cotangent bundle and is invariant under the
action of the geodesic flow (Theorem 2.1).

In general it will not be easy to find such Riemannian manifolds 
whose (punctured)cotangent bundle has a K\"ahler structure
where the symplectic form coincides with the
K\"ahler form and is invariant under the action of
the geodesic flow (see \cite{Sz1} and \cite{Sz2}).
Such a K\"ahler structure for complex and quaternion projective spaces
is constructed by making use of the Hopf fibration and the map above for
the sphere. 
Although the Cayley projective plane has no fiber bundle like the Hopf
fiber bundle, we prove here that a map similar to the cases of complex and
quaternion projective spaces gives an embedding of the punctured
cotangent bundle of the Cayley projective plane into the space of a
complexified exceptional Jordan algebra.
It is well-known that
the Cayley projective plane is one of the compact
symmetric spaces of rank one and that the exceptional Lie group
$F_4$ acts on it two-point homogeneously.  Some properties which we prove in
Theorem 2.1 could be shown quite easily if we used this property of minimal rank
for symmetric spaces, but we prove our main theorem through elementary
calculi in a Jordan algebra where the Cayley projective plane 
is realized as a subset consisting of primitive idempotents.
Further we give an embedding
of this image in the complexified exceptional Jordan algebra 
into the space of $8 \times 8$ complex matrices by composing with a map
given by \cite{Y}.  
 
In $\S 1$ we describe the Cayley projective plane as a subset
consisting of primitive idempotents in an exceptional Jordan algebra.
In $\S 2$ we state our main theorem and a corollary.
In $\S 3$ we recall some basic facts about the Jordan algebra 
including, so called, the Freudenthal product and the determinant 
on the Jodran algebra.
In $\S 4$  we prove our main theorem and in $\S 5$ we describe an
embedding of the punctured cotangent bundle of the Cayley projective
plane into the space of $8\times 8$ complex matrices.


\bigskip

\section{Cayley projective plane}
In this section we describe the Cayley projective plane as a subset
in the exceptional
Jordan algebra over the real number field $\mathbb{R}$(see \cite{M} and \cite{Be}).

Let $\mathbb{H}$ be the quaternion number field, that is, $\mathbb{H}$ is
an algebra over $\mathbb{R}$ generated by
$\{\bs{e}_i\}_{i=0}^3$ with the relations:
\begin{align}
\bs{e}_0\bs{e}_i &=\bs{e}_i\bs{e}_0\quad (i=0,1,2,3)\\
{\bs{e}_i}^2 &=-\bs{e}_0\quad (i=1,2,3)\\
\bs{e}_i\bs{e}_j &=-\bs{e}_j\bs{e}_i=\bs{e}_k \quad\pmod 3.
\end{align}

The Cayley number field $\mathcal{O}$ is a division algebra over $\mathbb{R}$
generated by $\{\bs{e}_i\}_{i=0}^7$ with $\bs{e}_i \bs{e}_j$ given by
the table
\begin{equation}
\begin{array}{|c|@{\kern0pt}|@{\kern0pt}|@{\kern0pt}|c|c|c|c|c|c|c|c|}
\hline
 & \bs{e}_0 & \bs{e}_1 & \bs{e}_2 & \bs{e}_3 & \bs{e}_4 & \bs{e}_5 &
\bs{e}_6 & \bs{e}_7\\
\HLINE
\bs{e}_0 & \bs{e}_0 & \bs{e}_1 & \bs{e}_2 & \bs{e}_3 & \bs{e}_4 &
\bs{e}_5 & \bs{e}_6 & \bs{e}_7\\
\hline
\bs{e}_1 & \bs{e}_1 & -\bs{e}_0 & \bs{e}_3 & -\bs{e}_2 & \bs{e}_5 &
-\bs{e}_4 & -\bs{e}_7 & \bs{e}_6\\
\hline
\bs{e}_2 & \bs{e}_2 & -\bs{e}_3 & -\bs{e}_0 & \bs{e}_1 & \bs{e}_6 &
\bs{e}_7 & -\bs{e}_4 & -\bs{e}_5\\
\hline
\bs{e}_3 & \bs{e}_3 & \bs{e}_2 & -\bs{e}_1 & -\bs{e}_0 & \bs{e}_7 &
-\bs{e}_6 & \bs{e}_5 & -\bs{e}_4\\
\hline
\bs{e}_4 & \bs{e}_4 & -\bs{e}_5 & -\bs{e}_6 & -\bs{e}_7 & -\bs{e}_0 &
\bs{e}_1 & \bs{e}_2 & \bs{e}_3\\
\hline
\bs{e}_5 & \bs{e}_5 & \bs{e}_4 & -\bs{e}_7 & \bs{e}_6 & -\bs{e}_1 &
-\bs{e}_0 & -\bs{e}_3 & \bs{e}_2\\
\hline
\bs{e}_6 & \bs{e}_6 & \bs{e}_7 & \bs{e}_4 & -\bs{e}_5 & -\bs{e}_2 &
\bs{e}_3 & -\bs{e}_0 & -\bs{e}_1\\
\hline
\bs{e}_7 & \bs{e}_7 & -\bs{e}_6 & \bs{e}_5 & \bs{e}_4 & -\bs{e}_3 &
-\bs{e}_2 & \bs{e}_1 & -\bs{e}_0\\
\hline
\end{array}
\end{equation}

Especially,
\begin{equation}
\text{$\bs{e}_1\bs{e}_4=\bs{e}_5,\,\bs{e}_2\bs{e}_4=\bs{e}_6$ and
$\bs{e}_3\bs{e}_4=\bs{e}_7$.}
\end{equation}
Hence $\mathcal{O}$ is identified with
\begin{equation}
\mathbb{H}\oplus\mathbb{H}\bs{e}_4
\end{equation}
and the multiplication between $x=a+b\bs{e}_4$ and $y=h+k\bs{e}_4$ $\in
\mathbb{H}\oplus\mathbb{H}\bs{e}_4$ is given by
\begin{equation}
x\cdot y=ah -\theta(k)b + \{ka +b\theta(h)\}\bs{e}_4,
\end{equation}
where $h=\sum_{i=0}^3 h_i\bs{e}_i\;(h_i\in\mathbb{R})$ and $\theta(h)
=h_0\bs{e}_0-h_1\bs{e}_1-h_2\bs{e}_2-h_3\bs{e}_3$, and so on.
We assume that the basis $\{\bs{e}_i\}_{i=0}^7$ are {\it {orthonormal}} 
{}and we will sometimes omit $\bs{e}_0 $ (= identity element)
and identify $\mathbb{R}
= \mathbb{R}\bs{e}_0 \subset \mathbb{H} \subset \mathcal{O}$.

For $h=\sum_{i=0}^7 h_i\bs{e}_i\in\mathcal{O}$, we denote
\begin{equation}
\theta(h)=h_0\bs{e}_0-\sum_{i=1}^7 h_i\bs{e}_i,
\end{equation}
as for $h \in \mathbb{H}$.

Let $M(3,\mathcal{O})$ be the space of $3\times 3$ matrices with
entries in $\mathcal{O}$. By identifying $M(3,\mathcal{O})\cong
M(3,\mathbb{R})\otimes_{\mathbb{R}}\mathcal{O}$, we denote for
$X\in M(3,\mathcal{O})$
\begin{equation}
\theta(X)=X_0\otimes\bs{e}_0-\sum_{i=1}^7 X_i\otimes \bs{e}_i,
\end{equation}
\begin{equation}
{}^t X=\sum_{i=0}^7{}^t X_i\otimes \bs{e}_i
\end{equation}
and
\begin{equation}
\tr X=\sum_{i=0}^7 (\tr X_i)\bs{e}_i,
\end{equation}
where $X=\sum_{i=0}^7X_i\otimes\bs{e}_i,\,X_i\in M(3,\mathbb{R})$.
\smallskip

Now consider the subspace $\mathfrak{J}$ in $M(3,\mathcal{O})$ 
defined by
\begin{equation}
\mathfrak{J}=\{X\in M(3,\mathcal{O})\;|\; \theta({}^t X)=X\},
\end{equation}
then $\dim_{\mathbb{R}}\mathfrak{J}=27$ and any $X\in\mathfrak{J}$ has
the form
\begin{equation*}
X=\left(\begin{array}{@{}ccc@{}}
\xi_1\bs{e}_0 & x_3 & \theta(x_2)\\
\theta(x_3) & \xi_2\bs{e}_0 & x_1\\
x_2 & \theta(x_1) & \xi_3\bs{e}_0
\end{array}\right)
\end{equation*}
where $\xi_i\in\mathbb{R}$, $x_i\in\mathcal{O}$.

The space $\mathfrak{J}$ is called an exceptional Jordan algebra with
the Jordan product
\begin{equation}
X\circ Y=\frac 12 (XY+YX)\in\mathfrak{J},
\end{equation}
$X,Y\in\mathfrak{J}$.

$\mathfrak{J}$ has an inner product given by
\begin{equation}
\tr (X\circ Y)=(X,Y)\bs{e}_0.
\end{equation}

In fact, for $X\in\mathfrak{J}$, we have $\tr X\in\mathbb{R}\bs{e}_0$, 
and 
\begin{equation}
(X,Y)= \sum_{i=1}^3(\xi_i\eta_i+2\left<x_i,y_i\right>)
\end{equation} 
where
\begin{align*}
X&=\left(\begin{array}{@{}ccc@{}}
\xi_1\bs{e}_0 & x_3 & \theta(x_2)\\
\theta(x_3) & \xi_2\bs{e}_0 & x_1\\
x_2 & \theta(x_1) & \xi_3\bs{e}_0
\end{array}\right)\\
Y&=\left(\begin{array}{@{}ccc@{}}
\eta_1\bs{e}_0 & y_3 & \theta(y_2)\\
\theta(y_3) & \eta_2\bs{e}_0 & y_1\\
y_2 & \theta(y_1) & \eta_3\bs{e}_0
\end{array}\right)
\end{align*}
and $\left<x_i,y_i\right>$ denote the inner product on $\mathcal{O}$.

Now the Cayley projective plane $P^2\mathcal{O}$ is defined as
\begin{defn}
$P^2\mathcal{O}=\{X\in\mathfrak{J}\;|\; X\circ X=X,\,\tr X= 1\}$.
\end{defn}

The exceptional Lie group $F_4$ is defined as a group of
algebra automorphisms of $\mathfrak{J}$, 
and acts on $P^2\mathcal{O}$ two-point homogeneously.
So we have
$P^2\mathcal{O}\cong F_4/\mathrm{Spin}(9)$, where $\mathrm{Spin}(9)$
is realized as a subgroup of $F_4$ consisting of those elements
which leave $\left(\begin{matrix}1 & 0 & 0\\
0 & 0& 0\\ 0& 0& 0\end{matrix}\right)\in\mathfrak{J}$ invariant.

\bigskip 

\section{A K\"ahler structure}

In this section we describe a K\"ahler structure on the punctured cotangent
bundle $T_0^* P^2\mathcal{O}$ = 
$T^*P^2\mathcal{O}\backslash P^2\mathcal{O}$ and state our main Theorem 2.1.

The tangent bundle $T P^2\mathcal{O}$ is identified with the subset in
$\mathfrak{J}\times\mathfrak{J}$ as
\begin{multline}
T P^2\mathcal{O}\\
=\left\{(X,Y)\in\mathfrak{J}\times\mathfrak{J}\;\left|\; X\circ X=X,
\tr X= 1, X\circ Y=\frac 12 Y  \right. \right\}.
\end{multline}

We introduce the Riemannian metric on $P^2\mathcal{O}$ such that
\begin{equation}
(Y_1,Y_2)_{\mathbb{P}} = \frac 12 \tr (Y_1\circ Y_2) = \frac 12
(Y_1,Y_2),
\end{equation}
where $(X,Y_1)$, $(X,Y_2)$ $\in$ $T P^2\mathcal {O}$.

Then from the inclusions $\mathbb{C}$ $\subset$ 
$\mathbb H$ $\subset$ $\mathcal {O}$, the complex projective plane $P^2\mathbb{C}$ and
the quaternion projective plane $P^2\mathbb{H}$ are
embedded isometrically into $P^2\mathcal{O}$ as totally geodesic 
submanifolds (\cite{Be}).

In the following we identify the tangent bundle $TP^2\mathcal{O}$ and
the cotangent bundle $T^*P^2\mathcal {O}$ through the metric
above. Under this identification the symplectic 
form $\omega_{\mathcal{O}}$ on $T^*P^2\mathcal{O}$
is given by
\begin{equation}
\omega_{\mathcal{O}}=-\frac 12(dX,dY),
\end{equation}
where we should interpret the inner product $(dX,dY)$ 
as a two-form in such a way that
\begin{equation}
-(dX,dY)
=\sum_{i=1}^3d\eta_i\wedge d\xi_i+2\sum_{i=1}^3\sum_{\alpha=0}^7
dy_{\alpha}^i\wedge dx_{\alpha}^i
\end{equation}
restricted to $TP^2\mathcal{O}$, that is, we notice that $\sum\eta_i\xi_i$
is replaced by $\sum d\eta_i\wedge d\xi_i$ and so on in the definition
of the inner product on $\mathfrak{J}$

We can extend $\theta:\mathfrak{J}\longrightarrow\mathfrak{J},\,
{}^t:\mathfrak{J}\longrightarrow\mathfrak{J},\tr:
\mathfrak{J}\longrightarrow\mathcal{O}$, and the inner
product $(\cdot,\cdot)$ to the
complexification $\mathfrak{J}\otimes_{\mathbb{R}}\mathbb{C} 
= \mathfrak{J}^{\mathbb{C}}$
in a natural way. So the Hermite inner 
product $\langle\cdot,\cdot\rangle$ on 
$\mathfrak{J}^{\mathbb{C}}$
is given by
\begin{equation}
\inner{X}{Y}
=(X,\overline{Y}),
\end{equation}
where $\overline{X}=\sum_{\alpha=0}^{7}\overline{X}_{\alpha}
\otimes\bs{e}_{\alpha}$, $X_{\alpha}\in M(3,\mathbb{C})$
and $\overline{X}_{\alpha}$ is the complex conjugate of
$X_{\alpha}$. The norm of these elements in $\mathfrak{J}$
and $\mathfrak{J}^\mathbb{C}$ is always written as $\norm \cdot$,
and we write the norm of the tangent vector $Y \in T_X(P^2\mathcal{O})$
by $\norm Y_{\mathbb{P}}$. 

Now consider the 
map $\tau_{\mathcal{O}}: T_0^* P^2\mathcal{O}(\cong T_0P^2\mathcal{O})
\longrightarrow \mathfrak{J}^\mathbb{C}$ defined by
\begin{multline}
\tau_{\mathcal{O}}(X,Y)\\
=(\norm{Y}^2X-Y\circ Y)\otimes 1+\frac{1}{\sqrt{2}}\norm{Y}Y\otimes\sqrt{-1}\\
= (2\norm{Y}^2_{\mathbb{P}}X-Y\circ Y)\otimes 1
+\norm{Y}_{\mathbb{P}}Y\otimes\sqrt{-1},
\end{multline}
then we have

\begin{thm}
The map $\tau_{\mathcal{O}}$ gives an isomorphism between
$T_0^*P^2\mathcal{O}$ and $\mathbb{E}
=\{A\in\mathfrak{J}^\mathbb{C}\;|\; A\circ A=0,A\ne 0\}$. Moreover
\begin{equation}
\tau^*_{\mathcal{O}}({\sqrt{-1}}\,\overline{\partial}\partial
\norm{A}^{\frac 12})=\frac{1}{\sqrt{2}}\omega_{\mathcal{O}}.
\end{equation}
\end{thm}

The two-form $\sqrt{-2}\,\overline{\partial}\partial
\norm{A}^{\frac 12}$ 
is itself a K\"ahler form on $\mathfrak{J}^{\mathbb{C}}\backslash\{0\}$, so
that we can regard $\mathfrak{J}^{\mathbb{C}}\backslash\{0\}$ is a
symplectic manifold. On this symplectic manifold
the flow $\{\phi_t\}_{t \in \mathbb{R}}$ 
defined by 
\begin{equation*}
\phi_t : A \mapsto \phi_t(A) = e^{-2{\sqrt {-1}}t}\cdot A
\end{equation*}
is a Hamilton flow. The Hamiltonian of this flow is given by
the function $f : A \mapsto \frac{1}{\sqrt 2}\norm{A}^{\frac 12}$.  
Since $\mathbb{E}$ is holomorphic
and the flow $\{\phi_t\}$ leaves $\mathbb{E}$ invariant,
the Hamiltonian of this flow on $\mathbb{E}$ is just
the restriction of $f$ to $\mathbb{E}$, that is, the Hamiltonian
is the square root of the metric function.  So the flow $\{\phi_t\}$ 
is the bicharacteristic flow of the square root of the Laplacian
on $P^2\mathcal {O}$.  Especially the flow
restricted to the unit sphere = $\{(X,Y) \in TP^2\mathcal {O}:
\norm{Y}_{\mathbb{P}}= 1\}$ coincides with the geodesic flow.  So we have
\begin{cor}
The geodesics $\gamma (t)$ on $P^2\mathcal{O}$ through a point $X$
with the direction $Y$ ($\norm{Y}_{\mathbb{P}}=1$ and 
$X\circ Y =\frac 12 Y$) is given by
\begin{equation} \gamma (t) = 
\cos 2t\cdot (X- \frac 12 Y\circ Y) + \frac {1}{2}\sin 2t\cdot Y +
\frac 12 Y\circ Y.
\end{equation}
\end{cor}

\bigskip
\section{Freudenthal product and determinant}

In this section we recall several formulas in the Jordan algebra
$\mathfrak{J}$ for later use(see \cite{M}).

Let $X,Y\in\mathfrak{J}$, then the ``Freudenthal product" $X\times Y$
$\in$ $\mathfrak{J}$
is defined by the formula
\begin{multline}
X\times Y=\frac 12 \left\{2X\circ Y-(\tr X)Y-(\tr Y)X\right.\\
\left. +(\tr X\cdot\tr Y-\tr(X\circ Y))E\right\}
\end{multline}
where $E=\left(\begin{matrix}1 & 0 & 0\\
0 & 1 & 0\\
0 & 0 & 1
\end{matrix}\right)$ 
and the determinant,``det $X$", for $X$ $\in$ $\mathfrak{J}$ 
is defined by
\begin{equation}
\det X=\frac 13 \tr (X\circ(X\times X)).
\end{equation}
Then we have
\begin{prop}
\begin{enumerate}
\renewcommand{\labelenumi}{(\roman{enumi})}    
\item $(X\circ Y,Z)= (X,Y\circ Z)$, for $\forall X,Y,Z$ $\in$ $\mathfrak{J}$ 
\item $X\circ (X\times X)=\det X\cdot E$(\text{Cayley-Hamilton})
\item $(X\times X)\times (X\times X)=\det X\cdot X$.
\end{enumerate}
\end{prop}

As we explained above
\begin{equation}
F_4=\{g\in\GL(\mathfrak{J})\;|\; g(X\circ Y)=g(X)\circ g(Y),\,
\forall X,Y\in\mathfrak{J}\}.
\end{equation}
Now $F_4$ is also given in the following ways:
\begin{align*}
F_4
&=\{g\in\GL(\mathfrak{J})\;|\;\text{for any $X,Y\in\mathfrak{J}$,
$\det(gX)=\det X,(gX,gY)=(X,Y)$}\}\\
&=\{g\in\GL(\mathfrak{J})\;|\;\text{for any $X,Y\in\mathfrak{J}$,
$\det(gX)=\det X, g(E)=E$}\}\\
&=\{g\in\GL(\mathfrak{J})\;|\;\text{for any $X,Y\in\mathfrak{J}$,
$g(X\times Y)=g(X)\times g(Y)$}\}.
\end{align*}
We also have for $g\in F_4$
\begin{equation}
\tr gX=\tr X.
\end{equation}

The ``Freudenthal product" and ``$\det$" on $\mathfrak{J}$ are
extended naturally to the complexification
$\mathfrak{J}^\mathbb{C}$, and we denote
them with the same notations. Then the complexification
of $F_4$ is defined in the same way:
\begin{defn}
The complex simple Lie group ${F_4}^{\mathbb{C}}$ is
\begin{align*}
{F_4}^{\mathbb{C}}
&=\{g\in\GL(\mathfrak{J}^\mathbb{C})
\;|\;g(X\circ Y)=g(X)\circ g(Y)\}\\
&=\{g\in\GL(\mathfrak{J}^\mathbb{C})
\;|\;\det(gX)=\det X,(g(X),g(Y))
=(X,Y)\}\\
&=\{g\in\GL(\mathfrak{J}^\mathbb{C})
\;|\;g(X\times Y)=g(X)\times g(Y)\}
\end{align*}
\end{defn}

The two-point  homogeneity of $F_4$ on $P^2\mathcal{O}$ is
equivalent to
\begin{prop}
Let $S(P^2\mathcal{O})=\{(X,Y)\;|\;(X,Y)\in TP^2\mathcal{O} \subset 
\mathfrak{J} \times \mathfrak{J},
\norm{Y}=1\}$, then $F_4$ acts on $S(P^2\mathcal{O})$
transitively and we have
\begin{equation}
S(P^2\mathcal{O})=F_4/\mathrm{Spin}(7),
\end{equation}
where the stationary subgroup at the point $\left(\left(
\begin{matrix}
1 & 0 & 0\\ 0 & 0 & 0\\0 & 0 & 0
\end{matrix}\right),\left(\begin{matrix}
0 & \frac{1}{\sqrt{2}} & 0\\ \frac{1}{\sqrt{2}} & 0 &0\\
0 & 0& 0\end{matrix}\right)\right)$ $\in$ 
$S(P^2\mathcal{O})$
is identified with
$\mathrm{Spin}(7)$.
\end{prop}

\begin{rem}
The representation of $F_4^{\mathbb{C}}$ on
$\mathfrak{J}^{\mathbb{C}}_0$(=$\{A \in \mathfrak{J}^{\mathbb{C}}| \tr
A = 0\}$)
is irreducible, and the subspace $\mathbb{E}$ is the orbit of 
the highest weight vector.
According to the theorem by
Lichtenstein(\cite{Li}), then such an orbit is characterized as the null
set of a certain system of quadric equations.  Thus our equation $A\circ A = 0$
(also this is equivalent to $A\times A =0$ and $\tr A = 0$)
is nothing else than an example of this theorem in \cite{Li}, however
we need the map $\tau_{\mathcal{O}}$ to put a K\"ahler structure on 
$T^*_0P^2\mathcal{O}$. 
\end{rem}

\bigskip
\section{Proof of the theorem}

We give a proof of Theorem 2.1 in a series of lemmas and
proposition.

First we prove
\begin{lem}
Let $(X,Y)\in TP^2\mathcal{O}$, then
\begin{enumerate}
\renewcommand{\labelenumi}{(\roman{enumi})}
\item $\det X=0$,
\item $\tr Y=0$,
\item $\det Y=0$.
\end{enumerate}
\end{lem}
\begin{proof}
Since
\begin{align*}
X\times X
=\frac 12\{2X\circ X-2\tr X\cdot X+((\tr X)^2-(X,X))E\}=0,
\end{align*} 
we have
\begin{equation*}
\det X = 0
\end{equation*}
by the definition of the determinant
\begin{equation*}
X\circ (X\times X)= \det X \cdot E = 0.
\end{equation*}

{}From the equality
\begin{equation*}
(X\circ Y,Z)=(X,Y\circ Z)\quad(X,Y,Z\in\mathfrak{J})
\end{equation*}
we have 
\begin{align*}
&(X\circ Y,X)=\frac 12 (Y,X)\\
=&\;(Y,X\circ X)=(Y,X).
\end{align*}
Hence $\tr Y=\dfrac 12 \tr(X\circ Y)=(X,Y)=0$.

Since $\tr Y = 0$, the Freudenthal product $Y\times Y$ is
expressed as
\begin{equation*}
Y\times Y=Y\circ Y-\frac 12 \norm{Y}^2 E,
\end{equation*}\\
and we have
\begin{equation}\label{e:freudeltal-Y}
Y\circ (Y\times Y)=Y\circ(Y\circ Y)-\frac 12 \norm{Y}^2 Y
=\det Y\cdot E.
\end{equation}
Then by taking the trace we have
\begin{equation*}
\tr(Y\circ (Y\circ Y))=3\det Y.
\end{equation*}
Now
\begin{align*}
&\tr(Y\circ (Y\circ Y))\\
=&\;\left(Y,Y\times Y+\frac 12 \norm{Y}^2 E\right )\\
=&\;2\left(X\circ Y,Y\times Y+\frac 12 \norm{Y}^2 E\right)\\
=&\;2\left(X,Y\circ\left(Y\times Y+\frac 12 \norm{Y}^2 E
\right)\right)\\
=&\;2 (X,Y\circ (Y\times Y))+\left( X, \norm{Y}^2 Y\right)\\
=&\;2 (X,\det Y\cdot E) = 2 \det Y.
\end{align*}
Hence we have
\begin{equation}
\det Y=0.
\end{equation}
\end{proof}

\begin{lem}
For $(X,Y) \in T_0P^2\mathcal O$, 
$\tau_{\mathcal O}(X,Y) \in \mathbb E$, 
that is, $\tau_{\mathcal O}(X,Y) \circ \tau_{\mathcal O}(X,Y) = 0$.
\end{lem}
\begin{proof}
Let $X,Y\in T_0P^2\mathcal O$, then
\begin{align*}
&\tau_{\mathcal{O}}(X,Y)\circ\tau_{\mathcal{O}}(X,Y)\\
=& \left((\norm{Y}^2X-Y\circ Y)^2-\frac 12 \norm{Y}^2 Y\circ
Y\right)\otimes 1\\
&\qquad +\frac{2}{\sqrt{2}}\norm{Y}Y\circ(\norm{Y}^2X-Y\circ Y)
\otimes\sqrt{-1}\\
=&\left(\norm{Y}^4X-2\norm{Y}^2 X\circ(Y\circ Y)+(Y\circ Y)\circ
(Y\circ Y)-\frac 12 \norm{Y}^2Y\circ Y\right)\otimes 1\\
&\qquad +\sqrt{2}\left(\frac 12 \norm{Y}^3 Y-\norm{Y}Y\circ(Y\circ
Y)\right)\otimes\sqrt{-1}.
\end{align*}

Here we notice the following formulas: let $(X,Y)\in
T(P^2\mathcal{O})$, then
\begin{enumerate}
\renewcommand{\labelenumi}{(\roman{enumi})}
\item $Y\circ(Y\circ Y)=\frac 12 \norm{Y}^2 Y$,
\item $(X + Y)\times(X + Y) = (X - Y)\times(X - Y)
=Y\circ Y-\frac 12\norm{Y}^2 E$,
\item  $\det (X\pm Y) = 0$,
\item  $X\circ (Y\circ Y)=\frac 12 \norm{Y}^2 X$.
\end{enumerate}
(i) is obtained by the Cayley-Hamilton
\begin{equation*}
Y\circ(Y\times Y)=\det Y \cdot E = 0.
\end{equation*}
and (ii) is easily shown.  (iii) and
(iv) are  proved by first calculating
\begin{align*}
&(X\pm Y)\circ((X\pm Y)\times(X\pm Y))\\
=&(X\pm Y)\circ\left(Y^2-\frac 12 \norm{Y}^2 E\right)\\
=&X\circ (Y\circ Y) - \frac 12 \norm{Y}^2 X\\
=&\det(X\pm Y) \cdot E,
\end{align*}
and by taking trace of both sides, we know $\det (X\pm Y) = 0$.  
Hence we have also (iv).

Next, from the formula $(Y\times Y)\times(Y\times Y)=(\det Y)Y=0$
we have
\begin{equation*}
0=\left(Y\circ Y-\frac 12 \norm{Y}^2E\right)\times
\left(Y\circ Y-\frac 12 \norm{Y}^2E\right)
\end{equation*}
and so we have
\begin{equation}
\norm{Y\circ Y}^2=\frac 12 \norm{Y}^4
\end{equation}
and
\begin{equation}
(Y\circ Y)\circ(Y\circ Y)=\frac 12 \norm{Y}^2 Y\circ Y.
\end{equation}

Finally by making use of these formulas we can prove
\begin{align*}
&=\tau_{\mathcal O}(X,Y) \circ \tau_{\mathcal O}(X,Y)\\
&=\left(\norm{Y}^4 X-2\norm{Y}^2\cdot\frac 12 \norm{Y}^2 X
+\frac 12 \norm{Y}^2 Y\circ Y-\frac 12 \norm{Y}^2 Y\circ Y\right)
\otimes 1\\
&\qquad+\sqrt{2}\left(\norm{Y}^3\cdot\frac 12 Y-\norm{Y}
\cdot\frac 12\norm{Y}^2 Y\right)\otimes\sqrt{-1}\\
&=0.
\end{align*}
\end{proof}

Let a map $\sigma:\mathbb {E}\longrightarrow \mathfrak{J}\times\mathfrak{J}$ be
\begin{equation}
\sigma: A \mapsto (X,Y),
\end{equation}
where $X$ and $Y$ are given by the following formulas:

\begin{align}
X
&=\frac 12 \frac{1}{\norm{A}}(A+\overline{A})+
\frac{A\circ \overline{A}}{\norm{A}^2},\\
Y
&=-\frac{\sqrt{-1}}{\sqrt{2}}\norm{A}^{-\frac 12}
(A-\overline{A}).
\end{align}
\begin{prop}
Let $X,Y$ be defined as above for $A\in\mathbb{E}$,
then
\begin{enumerate}
\renewcommand{\labelenumi}{(\roman{enumi})}
\item $X\circ X=X$, ${}^t\theta(X)=X$, $\tr X=1$
\item $Y\circ X=\frac 12 Y$, $\tr Y=0$.
\end{enumerate}
\end{prop}
\begin{proof}
By the definition of the Freudenthal product we have
\begin{align*}
A\times A
&=\frac 12\{2A\circ A-2\tr A\cdot A+((\tr A)^2-(A,A))E\}\\
&=\frac 12\{-2\tr A\cdot A+(\tr A)^2 E\},
\end{align*}
and
\begin{align*}
&A\circ (A\times A)=\det A \cdot E\\
=& A\circ\left(-\tr A\cdot A+\frac 12 (\tr A)^2 E\right)
=\frac 12 (\tr A)^2 A.
\end{align*}
So we have
\begin{equation}
\det A\cdot A=\frac 12 (\tr A)^2 A\circ A = 0
\end{equation}
and then we have 
\begin{equation}
\det A=0,\tr A=0.
\end{equation}

It follows easily that  $\,{}^t\theta(X)=X$
because of $\,{}^t\theta(A)=A$ and ${}^t\theta(\overline{A})=\overline{A}$.

Put $A=a\otimes 1+b\otimes\sqrt{-1}$ $\in$ $\mathbb{E}$, where $a,b\in\mathfrak{J}$.
Now we can assume $\norm{A}=1$, because of the homogeneity of the map $\sigma$.
Then we have  $a\circ a=b\circ b$, $a\circ b=0$, ${\norm a}^2 = {\norm b}^2
= \frac 12$,
and $\tr a = \tr b = 0 $.
Also we show
\begin{equation}
\det a=\det b = 0,\quad\det(a\pm b)=0.
\end{equation}
The last equalities are proved by the following argument:
from the equalities $a\times a=a\circ a-\frac 14 E$ 
and $b\times b=b\circ b-\frac 14 E$, 
we have
\begin{equation*}
a\circ (a\times a)=a\circ (a\circ a)-\frac 14 a=\det a\cdot E,
\end{equation*}
\begin{equation*}
b\circ (b\times b)=b\circ (b\circ b)-\frac 14 b=\det b\cdot E,
\end{equation*}
\begin{equation*}
(a\pm b)\times (a\pm b)=2 a\circ a-\frac 12 E=2b\circ b
-\frac 12 E.
\end{equation*}
Then 
\begin{align*}
(a\pm b)\circ ((a\pm b)\times(a\pm b))
&=(a\pm b)\circ\left(2 a\circ a -\frac 12 E\right)\\
&=(a\pm b)\circ\left(2 b\circ b -\frac 12 E\right)\\
&= 2a\circ (a\circ a -\frac 14 E)\pm 2b\circ (b\circ b -\frac 14 E) \\ 
&=\det(a\pm b)E.\\
\end{align*}
Hence we have
\begin{equation*}
2(\det a\pm\det b)=\det(a\pm b).
\end{equation*}
On the other hand
\begin{align}
&((a\pm b)\times (a\pm b))\times((a\pm b)\times (a\pm b))
\label{eq:x1}\\
=&\det(a\pm b)(a\pm b)\\
=&\left(2a\circ a-\frac 12 E\right)\times\left(2a\circ a
-\frac 12 E\right).
\end{align}
Note that the last equality shows that \eqref{eq:x1} does not
depend on the sign.
Hence
\begin{equation*}
\det(a+b)\cdot(a+b)=\det(a-b)\cdot (a-b).
\end{equation*}
So we have 
\begin{equation*}
\det(a+b)=\det(a-b),
\end{equation*}
since $\tr(a\circ a) = \frac 12$ and
\begin{equation}
\det a= \det b=0.
\end{equation}
Then these finally imply
\begin{equation}
\det(a\pm b)=0.
\end{equation}

By making use of these formulas
we prove $X\circ X=X$. Again we may assume $\norm{A}=1$, then
$X$ is written as
\begin{align*}
X &=a\otimes 1+(a\otimes 1+b\otimes \sqrt{-1})\circ
(a\otimes 1-b\otimes \sqrt{-1})\\
&=a+2a\circ a.
\end{align*}
Hence
\begin{align*}
X\circ X
&=(a+2a\circ a)\circ(a+2a\circ a)\\
&=a\circ a+4a\circ(a\circ a)+4(a\circ a)\circ(a\circ a)\\
&=a\circ a+4\cdot\frac a4+a\circ a=a+2a\circ a=X,
\end{align*}
where we used the equality
\begin{equation*}
(a\circ a)\circ(a\circ a)=\frac 14 a\circ a.
\end{equation*}

Also we have $\tr X=1=2\tr(a\circ a)$.

Next we show    
\begin{equation*}
X\circ Y=\frac 12 Y. 
\end{equation*}
Since
\begin{equation*}
Y=-\frac{\sqrt{-1}}{\sqrt{2}}(a\otimes 1+b\otimes\sqrt{-1}
-(a\otimes 1-b\otimes \sqrt{-1})) = \sqrt{2}b,
\end{equation*}
\begin{align*}
X\circ Y
&=\sqrt{2}(a+2a\circ a)\circ b\\
&=\sqrt{2}(a\circ b+2 b\circ (a\circ a))\\
&=2\sqrt{2}b\circ\left(a\times a+\frac 14 E\right)\\
&=2\sqrt{2}b\circ\left(b\times b+\frac 14 E\right)\\
&=2\sqrt{2}\left(\det b\cdot E+\frac 14 b\right)\\
&=\frac 12 \cdot(\sqrt{2}b)=\frac 12 Y.
\end{align*}
{}From these we have proved that $\tau_{\mathcal{O}}$ is a bijection
between $T_0^* P^2\mathcal{O}$ in $\mathfrak{J}\times\mathfrak{J}$
and $\mathbb{E}$ in $M(3,\mathcal{O})\otimes_{\mathbb{R}}
\mathbb{C}$ and that $\sigma$  is the inverse map.
\end{proof}

Next we prove $\tau_{\mathcal{O}}^*(\sqrt{-2}\,\overline{\partial}
\partial\sqrt{\norm{A}})=\omega_{\mathcal{O}}$, that is, the
K\"ahler form $\sqrt{-2}\,\overline{\partial}\partial\sqrt{\norm{A}}$
coincides with the symplectic form $\omega_{\mathcal{O}}$
on $T_0^*P^2\mathcal{O}$.

First we have
\begin{align*}
&\tau_{\mathcal{O}}^*(\sqrt{-1}\,\overline{\partial}
\partial\norm{A}^{\frac 12})\\
=&\tau_{\mathcal{O}}^*(\sqrt{-1}\,\overline{\partial}
\partial(A,\overline{A})^{\frac 14})\\
=&\frac{\sqrt{-1}}{4}d(\tau_{\mathcal{O}}^* 
(A,\overline{A})^{-\frac 34}
(dA,\overline{A}))\\
=&\frac{\sqrt{-1}}{4}d(\norm{Y}^{-3}(d\tau_{\mathcal{O}}^* A,
\tau_{\mathcal{O}}^*(\overline{A}))).
\end{align*}
Here we should consider $A\in\mathfrak{J}^{\mathbb{C}}$
to be the section
\begin{equation}
\begin{array}{ll}
A: & \mathbb{E}\longrightarrow
\mathbb{E}\times\mathfrak{J}^{\mathbb{C}}\\
& A\longmapsto (A,A)                      
\end{array}
\end{equation}
of the trivial bundle $\mathbb{E}\times\mathfrak{J}^{\mathbb{C}}$ 
on $\mathbb{E}$, and $dA$ the section of
$\mathfrak{J}^{\mathbb{C}}\otimes T^* \mathbb{E}$.
Note that the inner product $( \cdot, \cdot)$
defines the pairing $\mathfrak{J}\times \mathfrak{J}\otimes
T^*\mathbb{E}$ $\rightarrow$ $T^*\mathbb{E}$.  In the calculations
below we will 
use this pairing with the same notation $(\cdot,
\cdot)$.
Also we note $\norm{A}^2= \norm {Y}^4$ under the
mapping $\tau_{\mathcal{O}}$.

We write
\begin{align}
\tau_{\mathcal{O}}^*(A)
&=(\norm{Y}^2 X-Y\circ Y)\otimes 1+\frac{1}{\sqrt{2}}
\norm{Y}Y\otimes\sqrt{-1}\\
&=a\otimes 1+b\otimes\sqrt{-1},
\end{align}
where $a=a(X,Y),b=b(X,Y)$. Then
\begin{align*}
&\tau_{\mathcal{O}}^* (dA,\overline{A})\\
=& (d\tau_{\mathcal{O}}^* A,\tau_{\mathcal{O}}^*
\overline{A})\\
=&(da\otimes 1+d(b\otimes\sqrt{-1}),a\otimes 1-b\otimes
\sqrt{-1})\\
=&(da,a)
+(db,b)
+((db,a)
-(da,b))\sqrt{-1}.
\end{align*}
Now from
\begin{equation*}
(a,a)=\frac 12 \norm{Y}^4=(b,b)
\end{equation*}
we have
\begin{align*}
(da,a)&=(Y,Y)(dY,Y)\\
(db,b)&=(Y,Y)(dY,Y).
\end{align*}
So the real part of
\begin{equation*}
\tau_{\mathcal{O}}^* ((dA,\overline{A}))
\end{equation*}
is a closed form, since $d((Y,Y)(dY,Y))
=2(dY,Y)\wedge(dY,Y)=0$.
{}From
\begin{equation*}
(a,b)=\left(\norm{Y}^2 X-Y\circ Y,\frac{1}{\sqrt{2}}
\norm{Y}Y\right)=0,
\end{equation*}
\begin{align*}
&(db,a)-(b,da)\\
=&2(db,a)
=2\frac{1}{\sqrt{2}}\left(Y\otimes\frac{(dY,Y)}{\norm{Y}}
+\norm{Y}dY,\norm{Y}^2X-Y\circ Y\right)\\
=&\frac{2}{\sqrt{2}}\left\{\norm{Y}^3
(dY,X)-\norm{Y}(dY,Y\circ Y)
\right\}.
\end{align*}
Hence
\begin{align*}
&\frac{(\sqrt{-1})^2}{4}\cdot\frac{2}{\sqrt{2}}
\left[d\left\{\norm{Y}^{-3}\cdot (\norm{Y}^{3}(dY,X)
-\norm{Y}(dY,Y\circ Y))\right\}\right]\\
=&\frac{-1}{2\sqrt{2}}\left\{d(dY,X)
-d\left(\frac{(dY,Y\circ Y)}
{\norm{Y}^2}\right)\right\}.
\end{align*}
By the equality $(X,Y\circ Z) = (X\circ Y, Z)$ we have
\begin{equation*}
(dY,Y\circ Y)= (Y,d Y\circ Y)
\end{equation*}
and by $(Y,Y\circ Y) = 0$, we have
\begin{equation*}
(dY,Y\circ Y)=0.
\end{equation*}
Hence we finally proved
\begin{equation}
\tau_{\mathcal{O}}^* (\sqrt{-1}\,\overline{\partial}\partial
\sqrt{\norm{A}})=\frac{1}{2\sqrt{2}}(dY,dX)
= \frac{1}{\sqrt{2}}\omega_{\mathcal{O}}.
\end{equation}

\begin{rem}
The map $\tau_{\mathcal{O}}$ commutes with the actions of $F_4$ 
on $T_0P^2{\mathcal{O}}$ and on $\mathfrak{J}^{\mathbb{C}}$
(as a subgroup in $F_4^{\mathbb{C}})$.  Of course all the elements in
$F_4$ preserve the symplectic form $\omega_{\mathcal{O}}$.  Then
it will be true that the subgroup of $F_4^{\mathbb{C}}$ consisting
of those elements that preserve the symplectic form
$\omega_{\mathcal {O}}$ is compact. Hence it will coincide with $F_4$.
\end{rem}

\bigskip
\section{An embedding into $M(8,\mathbb{C})$}

In this section we describe an embedding of the space $\mathbb{E}
\subset\mathfrak{J}^\mathbb{C}$ into the space of $8\times 8$
complex matrices $M(8,\mathbb{C})$.

By identifying $\mathcal{O}\cong\mathbb{H}\oplus\mathbb{H}\bs{e}_4$,
we define $\gamma:\;\mathcal{O}\longrightarrow\mathcal{O}$ by
\begin{equation}
\gamma(h+k\bs{e}_4)=h-k\bs{e}_4.
\end{equation}
The map $\gamma$ is naturally extended to the complexification
$\mathcal{O}\otimes\mathbb{C}$, where we regard $\mathcal{O}\otimes
\mathbb{C}\cong \mathbb{H}\otimes\mathbb{C}\oplus\mathbb{H}\otimes
\mathbb{C}\bs{e}_4$. It is easily verified that $\gamma$ is an algebra
isomorphism of $\mathcal{O}$ (and of $\mathcal{O}\otimes\mathbb{C}$),
that is, $\gamma\in G_2$($\cong$ the group of algebra isomorphisms
of $\mathcal{O}$), $\gamma^2=\mathrm{Id}$, 
and $\theta\circ\gamma=\gamma\circ
\theta$.

Let $X\in\mathfrak{J}^\mathbb{C}$ be
\begin{equation*}
X=\left(\begin{array}{@{}ccc@{}}
\xi_1 & x_3 & \theta(x_2)\\
\theta(x_3) & \xi_2 & x_1\\
x_2 & \theta(x_1) & \xi_3
\end{array}\right),
\end{equation*}
$\xi_i\in\mathbb{C},x_i\in\mathcal{O}\otimes\mathbb{C},x_i=\sum_{\alpha =0}^7
x_i^{\alpha}\bs{e}_{\alpha}\,(x_i^{\alpha}\in\mathbb{C})$, and we
denote by
\begin{equation}
\gamma(X)=\left(\begin{array}{@{}ccc@{}}
\xi_1 & \gamma(x_3) & \theta(\gamma(x_2))\\
\theta(\gamma(x_3)) & \xi_2 & \gamma(x_1)\\
\gamma(x_2) & \theta(\gamma(x_1)) & \xi_3
\end{array}\right),
\end{equation}
then
\begin{equation}
\gamma:\;\mathfrak{J}^\mathbb{C}\longrightarrow\mathfrak{J}^\mathbb{C}
\end{equation}
is also an algebra isomorphism. 

We decompose elements $X\in \mathfrak{J}^\mathbb{C}$ as
\begin{align*}
X
&=\left(\begin{array}{@{}ccc@{}}
\xi_1 & x_3 & \theta(x_2)\\
\theta(x_3) & \xi_2 & x_1\\
x_2 & \theta(x_1) & \xi_3
\end{array}\right)\\
&=\left(\begin{array}{@{}ccc@{}}
\xi_1 & m_3 & \theta(m_2)\\
\theta(m_3) & \xi_2 & m_1\\
m_2 & \theta(m_1) & \xi_3
\end{array}\right)+
\left(\begin{array}{@{}ccc@{}}
0 & a_3 & -a_2\\
-a_3 & 0 & a_1\\
a_2 & -a_1 & 0
\end{array}\right)\bs{e}_4\\
&= M+A\bs{e}_4
\end{align*}
where $x_i=m_i+a_i\bs{e}_4\,(m_i,a_i\in\mathbb{H}\otimes\mathbb{C})$.
Then we can regard
\begin{equation*}
\mathfrak{J}\otimes\mathbb{C}=\mathfrak{J}(3,\mathbb{H})\otimes\mathbb{C}
\oplus\mathbb{H}^3\otimes\mathbb{C}.
\end{equation*}
Here we denote by $\mathfrak{J}(3,\mathbb{H})$ the Jordan algebra of
$3\times 3$ Hermitian matrices with entries in $\mathbb{H}$.

Now we define a map $g$ following \cite{Y}:
\begin{equation}
\begin{CD}
g:\; \mathfrak{J}\otimes\mathbb{C} @>>> \mathfrak{J}(4,\mathbb{H})_0
\otimes \mathbb{C}\\
\phantom{g:\;}\varin @. \varin\\
X=M+A\bs{e}_4 \;|\kern-.74em @>>> 
\left(\begin{array}{@{}cc@{}}
\frac 12 \tr M & \begin{array}{@{}ccc@{}} ia_1 & ia_2 & ia_3\end{array}\\
\begin{array}{@{}c@{}}i\theta(a_1)\\ i\theta(a_2)\\ i\theta(a_3)\end{array}
& M-\frac 12 \tr M\cdot E
\end{array}\right)
\end{CD}
\end{equation}
where we denote by $\mathfrak{J}(4,\mathbb{H})_0$ the subspace in the
Jordan algebra $\mathfrak{J}(4,\mathbb{H})$ consisting of 
elements whose trace is zero.

The map $g$ satisfies
\begin{enumerate}
\renewcommand{\labelenumi}{(\roman{enumi})}
\item $g(X)\circ g(Y)
= g(\gamma(X\times Y))+\frac 14(\gamma(X),Y)\cdot E$
\item $(g(X),g(Y))=(\gamma(X),Y)=\tr(\gamma(X)\circ Y)$

\noindent $=\tr(g(X)\circ g(Y))$
\end{enumerate}

We remark that $\tr (A\circ B)$ for $A,B\in\mathfrak{J}(4,\mathbb{H})$
defines a Euclidean inner product on $\mathfrak{J}(4,\mathbb{H})$.  We 
extend it to the complexification $\mathfrak{J}(4,\mathbb{H})\otimes
\mathbb{C}$.

If $A\in\mathfrak{J}^{\mathbb{C}}, A\circ A=0$, then from (i) above
we have at once

\begin{equation}
g(A)\circ g(A)=\frac 14(\gamma(A),A)\cdot E.
\end{equation}


Let $\rho:\;\mathbb{H}\otimes \mathbb{C}\rightarrow M(2,\mathbb{C})$ be the
isomorphism given by
\begin{equation}
\rho\left(\sum_{i=0}^3 z_i\bs{e}_i\right)\longmapsto
\left(\begin{array}{@{}cc@{}}
z_0+z_1i & z_2+z_3 i\\
-z_2+z_3i & z_0-z_1i
\end{array}\right),
\end{equation}
and we denote with the same notation $\rho$ the map
\begin{equation}
\begin{CD}
\rho:\;\mathfrak{J}(4,\mathbb{H})\otimes \mathbb{C} @>>> M(8,\mathbb{C})\\
\phantom{\rho:\;\mathfrak{J}\!}\varin @. \varin\\
\phantom{\rho:\;\mathfrak{J}}(h_{ij})\phantom{\otimes\mathbb{C}}
\;|\kern-1.9em @>>> (\rho(h_{ij}))
\end{CD}
\end{equation}
where $h_{ij}\in\mathbb{H}\otimes\mathbb{C}$.

\begin{prop}
\begin{equation}
\rho(\mathfrak{J}(4,\mathbb{H}))=\{A\in M(8,\mathbb{C})\;|\;
\mathbb{J}A={}^tA\mathbb{J}\},
\end{equation}
where
\begin{equation*}
\mathbb{J}=\left(\begin{array}{@{}cccc@{}}
J & 0 & 0 & 0\\
0 & J & 0 & 0\\
0 & 0 & J & 0\\
0 & 0 & 0 & J
\end{array}\right),\quad J=\left(\begin{array}{@{}cc@{}}
0 & 1\\ -1 & 0
\end{array}\right).
\end{equation*}
\end{prop}

Let $A\in\mathfrak{J}\otimes\mathbb{C}$ and $A\circ A=0$, then as we
know $\rho(g(A))^2$ is a scalar matrix and $\tr\rho(g(A))=0$.
Conversely we have

\begin{prop}
Let $A\in M(8,\mathbb{C}),\, A\ne 0,\,\mathbb{J}A={}^t A\mathbb{J}$, and
$A^2=\lambda E,\,\lambda\in\mathbb{C},\,\tr A=0$, then there exists
an element $X\in\mathfrak{J}^\mathbb{C}$ such that
\begin{align*}
X\circ X &= 0\\
\rho(g(X)) &= A
\end{align*}
under the condition that $A$ is of the form
\begin{equation}
A=\left(\begin{array}{@{}c|c|c|c@{}}
\begin{array}{cc} 0 & 0 \\ 0 & 0
\end{array}
&&&\\
\hline
 & \begin{array}{cc} \xi_1 & 0 \\ 0 & \xi_1\end{array} &&\\
\hline
&& \begin{array}{cc} \xi_2 & 0 \\ 0 & \xi_2\end{array} &\\
\hline
&&& \begin{array}{cc} \xi_3 & 0 \\ 0 & \xi_3\end{array}
\end{array}\right),\quad \sum_{i=1}^3 \xi_i=0.
\end{equation}
\end{prop}

\begin{rem}
The canonical line bundle of this complex structure on
$T^*_0P^2\mathcal{O}$ will be holomorphically trivial 
and this realization of $T^*_0P^2\mathcal{O}$ in the space
$M(8,\mathbb{C})$ given in the above proposition 
will be useful to construct an explicit holomorphic trivialization of 
the canonical line bundle.
\end{rem}

\bigskip

\providecommand{\bysame}{\leavevmode\hbox to3em{\hrulefill}\thinspace}


\begin{thebibliography}{Raw79}

\bibitem[Be]{Be}
A.~L. Besse, Manifolds all of whose Geodesics are closed,
Springer-Verlag, 1978.

\bibitem[FT]{FT}
K.~Furutani and R.~Tanaka, \textit{A K\"{a}hler structure on the punctured
cotangent bundle of complex and quaternion projective spaces and its
application to geometric quantization I}, J. Math. Kyoto Univ.
\textbf{34-4} (1994), 57--69.

\bibitem[FY]{FY}
K.~Furutani and S.~Yoshizawa, \textit{A K\"{a}hler structure on the punctured
cotangent bundle of complex and quaternion projective spaces and its
application to geometric quantization II}, Japanese J. Math. \textbf{21}
(1995), 355--392.

\bibitem[Ii1]{Ii1}
K.~Ii, \textit{On a Bargmann-type transform and a Hilbert space of
holomorphic functions}, T\^{o}hoku Math. J. (1) \textbf{38} (1986), 57--69.

\bibitem[Ii2]{Ii2}
K.~Ii and T.~Morikawa, \textit{K\"{a}hler structures on tangent bundle of Riemannian manifolds
of constant positive curvature}, Bull. Yamagata Univ. Natur. Sci.,
\textbf{14}(1999), no. 3, 141--154.

\bibitem[Li]{Li}
W.~Lichtenstein, \textit{A system of quadrics describing the orbit of the
highest weight vector}, Proceedings of Amer. Math. Soc., \textbf{84}(1982),
No. 4, 605--608.

\bibitem[M]{M}
S.~Murakami, \textit{Exceptional Simple Lie Groups and Related Topics
in Recent Differential Geometry, in Differential Geometry and
Topology}, Sprinter Lecture Note, vol. 1369, Proceedings, at Tianjir,
1986-1987.

\bibitem[Ra1]{Ra1}
J.~H. Rawnsley, \textit{Coherent states and K{\"a}hler manifolds}, Quart. J.
Math. Oxford \textbf{28} (1977), 403--415.

\bibitem[Ra2]{Ra2}
\bysame, \textit{A non-unitary pairing of polarization for the Kepler
problem}, Trans. Amer. Math. Soc. \textbf{250} (1979), 167--180.

\bibitem[So]{So}
J.~M. Souriau, \textit{Sur la vari\'{e}t\'{e} de Kepler}, Symposia Math.
\textbf{14} (1974), 343--360.

\bibitem[Sz1]{Sz1}
R.~Sz\H{o}ke, \textit{Complex structures on the tangent bundle of Riemannian
manifolds}, Math. Ann. \textbf{291} (1991), 409--428.

\bibitem[Sz2]{Sz2}
\bysame, \textit{Adapted Complex Structures and Geometric
Quantization}, Nagoya J. Math. \textbf{154} (1999), 171--183.

\bibitem[Y]{Y}
I.~Yokota, \textit{Realizations of involutive automorphisms $\sigma$ and
  $G^{\sigma}$ of exceptional linear Lie Groups $G$. Part $I$, $G =
  G_2$, $F_4$ and
 $ E_6$}, Tsukuba J. Math. \textbf{14} (1990), no. 1, 185--223.
\end{thebibliography}
\end{document}